\input amstex
\documentstyle{amsppt}
\magnification=\magstep1                        
\hsize6.5truein\vsize8.9truein                  
\NoRunningHeads
\loadeusm

\magnification=\magstep1                        
\hsize6.5truein\vsize8.9truein                  
\NoRunningHeads
\loadeusm

\document

\topmatter

\title
On the oscillation of the modulus of the Rudin-Shapiro polynomials around the middle of their ranges
\endtitle

\rightheadtext{On the oscillation of the modulus of the Rudin-Shapiro polynomials around the middle of their ranges} 

\author Tam\'as Erd\'elyi
\endauthor

\address Department of Mathematics, Texas A\&M University,
College Station, Texas 77843, College Station, Texas 77843 
\endaddress

\endthanks


\date
September 29, 2025  
\enddate
 
\email terdelyi\@tamu.edu
\endemail

\abstract
Let either $R_k(t) := |P_k(e^{it})|^2$ or $R_k(t) := |Q_k(e^{it})|^2$, where $P_k$ and $Q_k$ are 
the usual Rudin-Shapiro polynomials of degree $n-1$ with $n=2^k$. The graphs of the trigonometric polynomials 
$R_k$ on the period suggest many zeros of $R_k(t)-n$ in a dense fashion on the period. Let ${\Cal N}(I,R_k-n)$ 
denote the number of zeros, counted with multiplicities, of the trigonometric polynomial $R_k-n$ in an interval 
$I := [\alpha,\beta] \subset [0,2\pi)$. Improving earlier results proved only for the interval $I := [0,2\pi)$, 
in this paper we show that 
$$\frac{n|I|}{8\pi} - \frac{2}{\pi} (2n\log n)^{1/2} - 1 \leq N(I,R_k-n) \leq \frac{n|I|}{\pi} +  \frac{8}{\pi}(2n\log n)^{1/2}\,,\qquad k \geq 2\,,$$
for every interval $I := [\alpha,\beta] \subset [0,2\pi)$, where $|I| = \beta-\alpha$ denotes the length of the 
interval $I$.   
\endabstract

\endtopmatter 

\head 1. Introduction \endhead
Let $D := \{z \in {\Bbb C}: |z| < 1\}$ denote the open unit disk of the complex plane.
Let $\partial D :=  \{z \in {\Bbb C}: |z| = 1\}$ denote the unit circle of the complex plane. 
Littlewood polynomials are polynomials with each of their coefficients in $\{-1,1\}$. A special 
sequence of Littlewood polynomials are the Rudin-Shapiro polynomials, They appear in Harold 
Shapiro's 1951 thesis [17] at MIT and are sometimes called just the Shapiro polynomials. 
They also arise independently in Golay's paper [14]. They are remarkably simple to construct 
and are a rich source of counterexamples to possible conjectures. The Rudin-Shapiro polynomials 
are defined recursively as follows:
$$\split P_0(z) & :=1\,, \qquad Q_0(z) := 1\,, \cr 
P_{k+1}(z) & := P_k(z) + z^{2^k}Q_k(z)\,, \cr
Q_{k+1}(z) & := P_k(z) - z^{2^k}Q_k(z)\,, \cr \endsplit$$
for $k=0,1,2,\ldots\,.$ Note that both $P_k$ and $Q_k$ are polynomials of degree $n-1$ with $n := 2^k$
having each of their coefficients in the set $\{-1,1\}$.
It is well known and easy to check by using the parallelogram law that
$$|P_{k+1}(z)|^2 + |Q_{k+1}(z)|^2 = 2(|P_k(z)|^2 + |Q_k(z)|^2)\,, \qquad z \in \partial D\,.$$
Hence
$$|P_k(z)|^2 + |Q_k(z)|^2 = 2^{k+1} = 2n\,, \qquad z \in \partial D\,. \tag 1.1$$
It is also well known (see Section 4 of [3], for instance), that
$$Q_k(-z) = (-1)^{k+1}P_k^*(z) := (-1)^{k+1} z^{n-1}P_k(1/z)\,, \qquad k \geq 1\,, \quad z \in {\Bbb C} \setminus \{0\}\,,$$
and hence
$$|Q_k(-z)| = |P_k(z)|\,, \qquad z \in \partial D\,. \tag 1.2$$
Various properties of the Rudin-Shapiro polynomials are discussed in [4] and [5]. 
As for $k \geq 1$ both $P_k$ and $Q_k$ have odd degree,
both $P_k$ and $Q_k$ have at least one real zero. The fact that for $k \geq 1$ both $P_k$ and $Q_k$
have exactly one real zero was proved in [4]. It has been shown in [8] that
the Mahler measure (geometric mean) and the maximum modulus of the Rudin-Shapiro polynomials
$P_k$ and $Q_k$ of degree $n-1$ with $n := 2^k$ on the unit circle of the complex plane
have the same size. That is, in addition to (1.1), the Mahler measure of the Rudin-Shapiro polynomials
of degree $n-1$ with $n := 2^k$ is bounded from below by $cn^{1/2}$, where $c > 0$ is
an absolute constant. In [9] various results on the zeros of the Rudin-Shapiro polynomials 
are proved and some open problems are raised. In [10] a conjecture of Saffari on the asymptotic value of
the Mahler measure of the Rudin-Shapiro polynomials $P_k$ and $Q_k$ is proved to be $(2n/e)^{1/2} = (2^{k+1}/e)^{1/2}$.

For a monic polynomial 
$$P(z) = \prod_{j=1}^n{(z-\alpha_j)} = z^n + \sum_{j=0}^{n-1}{a_jz^j}\,, \qquad a_j \in {\Bbb C}\,, \quad a_0 \neq 0\,, \tag 1.3$$ 
let $$H(P) := \frac{1}{|a_0|^{1/2}} \max_{z \in \partial D}{|P(z)|}\,.$$  
Let  
$$\alpha_j = \rho_j e^{i\theta_j}\,, \qquad \rho_j > 0\,, \quad \theta_j \in [0,2\pi)\, \quad j \in \{1,2,\ldots,n\}\,. \tag 1.4$$
By using the notation (1.4), for a polynomial $P$ of the form (1.3) and the interval $I := [\alpha,\beta] \subset [0,2\pi)$ 
let $N(I,P)$ denote the number of the values $j \in \{1,2,\ldots,n\}$ for which $\theta_j \in I$. In 1950 Erd\H os and Tur\'an [14] 
proved the following result.

\proclaim{Theorem 1.1} We have 
$$\left| N(I,P) - \frac{n|I|}{2\pi} \right| \leq 16 (n\log H(P))^{1/2}$$
for every monic polynomial of the form (1.3) and for every interval $I := [\alpha,\beta] \subset [0,2\pi)$, where 
$|I| = \beta - \alpha$ denotes the length of the interval $I$.
\endproclaim

In [18] K. Soundararajan proved that the constant $16$ in the above result may be replaced by $8/\pi$. 
Moreover, using the notation $\log^+(x) := \max\{\log x,0\}$ for $x > 0$ he showed that the upper bound in Theorem 1.1 can be 
replaced by $(8/\pi)(nh(P))^{1/2}$, where
$$h(p) := \frac{1}{2\pi} \int_{0}^{2\pi}{\log^+ \frac{P(e^{it})}{|a_0]^{1/2}} \,dt} \leq H(P)\,.$$
Rudin-Shapiro polynomials play a key role in [2] as well as in [12] to prove the existence of flat Littlewood polynomials, 
a recent breakthrough result. More on Rudin-Shapiro polynomials may be found in [6,7,16].

\head 2. New Results \endhead
Let either $R_k(t) := |P_k(e^{it})|^2$ or $R_k(t) := |Q_k(e^{it})|^2$, and $n := 2^k$. 
In [1] we combined close to sharp upper bounds for the modulus of the autocorrelation coefficients of 
the Rudin-Shapiro polynomials with a deep theorem of Littlewood (see Theorem 1 in [15]) to prove that
there is an absolute constant $c>0$ such that the equation $R_k(t) = (1+\eta)n$ with $n := 2^k$ has at least 
$cn^{0.5394282}$ distinct solutions in $[0,2\pi)$ whenever $\eta$ is real, $|\eta| \leq  2^{-8}$, and $n$ is 
sufficiently large. In this paper we improve this result substantially. Let ${\Cal N}(I,R_k-n)$ denote the number
of zeros, counted with multiplicities, of the trigonometric polynomial $R_k(t)-n$ in an interval 
$I := [\alpha,\beta] \subset [0,2\pi)$.    

\proclaim{Theorem 2.1} Let $k \geq 0$ and $n := 2^k$ be integers. We have 
$$\frac{n|I|}{8\pi} - \frac{2}{\pi}(2n\log n)^{1/2} - 1 \leq {\Cal N}(I,R_k-n) 
\leq \frac{n|I|}{\pi} + \frac{8}{\pi}(2n\log n)^{1/2}\,,\qquad k \geq 2\,,$$
for every interval $I := [\alpha,\beta] \subset [0,2\pi)$, where $|I| = \beta-\alpha$ denotes the length of the interval $I$.  
\endproclaim

This extends the main result in [12] from the case of the interval $I := [0,2\pi)$ to the case of the interval 
$I = [\alpha,\beta] \subset [0,2\pi]$. 
In our proof of Theorem 2.1 we combine ideas used in [11] and a classical result of Erd\H os and Tur\'an 
[13] with a constant improved recently by Soundararajan [19]. 

\head 3. Lemmas \endhead

In the proof of Theorem 2.1 we need the lemma below stated and proved as Lemma 3.1 in [9].   

\proclaim{Lemma 3.1} Let $k \geq 2$ and $n := 2^k$ be integers, and let
$$z_j := e^{it_j}\,, \quad t_j := \frac{2\pi j}{n}\,, \quad j \in {\Bbb Z}\,.$$
We have
$$\split P_k(z_j) & = 2P_{k-2}(z_j)\,, \qquad j=2u\,, \enskip u \in {\Bbb Z}\,, \cr
P_k(z_j) & =(-1)^{(j-1)/2} 2i\,Q_{k-2}(z_j)\,, \qquad j=2u+1\,, \enskip u \in {\Bbb Z}\,, \cr \endsplit$$
where $i$ is the imaginary unit.
\endproclaim

For a trigonometric polynomial $T$ of the form
$$T(\theta) = \pm 2 \cos(m\theta) + \sum_{j=-m+1}^{m-1}{a_je^{ij\theta}}\,, \qquad a_j \in {\Bbb C}\,, \tag 3.1$$
let
$$H(T) := \max_{\theta \in \Bbb R}{|T(\theta)|}\,.$$
For an interval $I := [\alpha,\beta] \subset [0,2\pi)$ and a trigonometric polynomials $T$ of the form (3.1) let ${\Cal N}(I,T)$ 
denote the number of zeros, counted with multiplicities, of $T$ in $I$. 

\proclaim{Lemma 3.2}
We have
$${\Cal N}(I,T) - \frac{m|I|}{\pi} \leq \frac{8}{\pi} (2m \log H(T))^{1/2}$$
for every trigonometric polynomial $T$ of the form (3.1) and for every interval $I := [\alpha,\beta] \subset [0,2\pi]$, 
where $|I|:= \beta-\alpha$. \endproclaim

\demo{Proof}
This follows from the Erd\H os-Tur\'an inequality (Theorem 1.1) with $16$ replaced by Soundararajan's constant $8/\pi$. 
\qed \enddemo

\proclaim{Lemma 3.3} Let $k \geq 0$ and $n:=2^k$ be integers. We have
$${\Cal N}(I,R_k-n) - \frac{n|I|}{\pi} \leq \frac{8}{\pi}(2n\log n)^{1/2}$$
for every interval $I := [\alpha,\beta] \subset [0,2\pi)$, where $|I|:= \beta-\alpha$.
\endproclaim

\demo{Proof}
Observe that $R_k-n$ is of the form (3.1) with $m:=n-1$. It follows from (1.1) that
$$H(R_k-n) = \max_{\theta \in {\Bbb R}}|R_k(\theta)-n| \leq n\,,$$
and the lemma follows from Lemma 3.2
\qed \enddemo

Replacing $n$ by $n/4$ we get the following corollary.

\proclaim{Lemma 3.4} 
Let $k \geq 2$ and $n:=2^k$ be integers. We have
$$N(I,R_{k-2}-n/4) - \frac{n|I|}{4\pi} \leq \frac{4}{\pi}(2n\log n)^{1/2}$$
for every interval $I := [\alpha,\beta] \subset [0,2\pi]$, where $|I|:= \beta-\alpha$.
\endproclaim

\head 4. Proof of Theorem 2.1 \endhead

\demo{Proof of Theorem 2.1}
Let $k \geq 2$ and $n:=2^k$ be integers, and let $I := [\alpha,\beta] \subset [0,2\pi)$. 
Assume that $R_k(t) = |P_k(e^{it})|^2$. The case $R_k(t) = |Q_k(e^{it})|^2$ follows 
from it by (1.2). The upper bound of the theorem follows from Lemma 3.3.
We now prove the lower bound of the theorem, which is more subtle. Without loss of generality 
we may assume that 
$$|I| \geq \frac{4\pi}{n}\,,$$
otherwise the lower bound of the theorem is trivial. For the sake of brevity let
$$A_j := R_{k-2}(t_j)-n/4\,, \qquad j \in {\Bbb Z}\,,$$  
where $t_j:=2\pi j/n$ is the same as in Lemma 3.1. We define the integers $h$ and $M$ by   
$$t_h < \alpha \leq t_{h+1} < t_{h+M+1} \leq \beta < t_{h+M+2}\,.$$
Observe that 
$$M \geq \frac{n|I|}{2\pi} - 2\,. \tag 4.1$$
We study the $M$-tuple $\langle A_{h+1},A_{h+2},\ldots,A_{h+M} \rangle$.   
Lemma 3.4 implies that $R_{k-2}(t)-n/4$ has at most 
$$\frac{n|I|}{4\pi} + \frac{4}{\pi} (2n \log n)^{1/2} \tag 4.2$$ 
zeros in $I$. Therefore the Intermediate Value Theorem yields that the number of sign changes in the $M$-tuple 
$\langle A_{h+1},A_{h+2},\ldots,A_{h+M} \rangle$ is at most as large as the value in (4.2). 
Hence (4.1) and (4.2) imply that there are integers 
$$h+1 \leq j_1 < j_2 < \cdots < j_N \leq h+M$$ 
with 
$$N \geq \frac{n|I|}{2\pi} - 2 - \frac{n|I|}{4\pi} - \frac{4}{\pi}(2n\log n)^{1/2}
=\frac{n|I|}{4\pi} - \frac{4}{\pi} (2n \log n)^{1/2} - 2\, \tag 4.3$$ 
such that 
$$A_{j_{\nu}}A_{j_{\nu}+1} \geq 0\,, \qquad \nu = 1,2,\ldots,N\,. \tag 4.4$$   
Using Lemma 3.1 we have either
$$\split 16A_{j_{\nu}}A_{j_{\nu}+1} = & (4(R_{k-2}(t_{j_{\nu}})-n/4))(4(R_{k-2}(t_{j_{\nu}+1})-n/4)) \cr  
= & (4|P_{k-2}(e^{it_{j_{\nu}}})|^2-n)(4|P_{k-2}(e^{it_{j_{\nu}+1}})|^2-n) \cr 
= & (|P_k(e^{it_{j_{\nu}}})|^2-n)(|Q_k(e^{it_{j_{\nu}+1}})|^2-n) \cr
= & (|P_k(e^{it_{j_{\nu}}})|^2-n)(n-|P_k(e^{it_{j_{\nu}+1}})|^2)\,, \endsplit \tag 4.5$$
or
$$\split 16A_{j_{\nu}}A_{j_{\nu}+1} = & (4(R_{k-2}(t_{j_{\nu}})-n/4))(4(R_{k-2}(t_{j_{\nu}+1})-n/4)) \cr  
= & (4|P_{k-2}(e^{it_{j_{\nu}}})|^2-n)(4|P_{k-2}(e^{it_{j_{\nu}+1}})|^2-n) \cr 
= & (|Q_k(e^{it_{j_{\nu}}})|^2-n)(|P_k(e^{it_{j_{\nu}+1}})|^2-n) \cr 
= & (n-|P_k(e^{it_{j_{\nu}}})|^2)(|P_k(e^{it_{j_{\nu}+1}})|^2-n)\,. \endsplit \tag 4.6$$
Combining (4.4), (4.5), and (4.6), we can deduce that
$$(|P_k(e^{it_{j_{\nu}}})|^2-n)(|P_k(e^{it_{j_{\nu}+1}})|^2-n) = -16A_{j_{\nu}}A_{j_{\nu}+1} \leq 0\,, 
\qquad \nu = 1,2,\ldots,N\,.$$ 
Hence the Intermediate Value Theorem implies that $R_k(t)-n = |P_k(e^{it})|^2-n$ has at least one zero in 
each of the intervals 
$$[t_{j_{\nu}},t_{j_{\nu}+1}]\,, \qquad \nu = 1,2,\ldots,N\,.$$
Recalling (4.3) we conclude that $R_k(t)-n = |P_k(e^{it})|^2-n$ has at least 
$$N/2 \geq \frac{n|I|}{8\pi} - \frac{2}{\pi} (2n \log n)^{1/2}-1$$ 
distinct zeros in $I$.   
\qed \enddemo

\enddocument